\documentclass[12pt,oneside]{article}
\usepackage{mathtext}
\usepackage[T2A]{fontenc}
\usepackage[cp1251]{inputenc}
\usepackage[russian, english]{babel}
\usepackage[dvips]{graphicx,rotating}
\usepackage{amssymb}
\usepackage{amsthm}
\usepackage{bm}
\usepackage{latexsym}
\usepackage{titlesec}
\usepackage{amsmath}
\usepackage{amsfonts}
\usepackage{cite}
\usepackage{eufrak}
\usepackage{eucal}
\usepackage{color}
\usepackage[a4paper, left=20mm, right=20mm, bottom=20mm, head=20mm]{geometry}
\newcounter{example}[section]
\newcommand{\exam}{\addtocounter{example}{1}}

\begin{document}

\begin{center}
\textbf{\Large{Multiplicative inverse  for Wolff-Denjoy Series}}
\end{center}

\begin{center}
\textbf{ A. R. Mirotin, A. A. Atvinovskii}\\
amirotin@yandex.ru, aatvinovskiy@gmail.com
\end{center}

Abstract. Let a function $f$ with real poles that form a monotonic and bounded sequence be expanded in a Wolff-Denjoy series with positive coefficients. The main result of the note states that if we subtract its “linear part” from the function $1/f$, then the remaining “fractional part”\  of this function will also expand into Wolff-Denjoy series (its poles are also real, and the coefficients of the series are negative).
 Application of the result to operator theory is given.


\begin{flushleft}
УДК 517.9
\end{flushleft}

\begin{center}
\textbf{\large{О МУЛЬТИПЛИКАТИВНОМ ОБРАЩЕНИИ\footnote{От английского "multiplicative inverse".  Авторы отдают себе отчет в том, что термин "мультипликативное обращение функции $f$" применительно к функции $1/f$ не принят в математической литературе на русском языке.} РЯДОВ ВОЛЬФА-ДАНЖУА}}
\end{center}

\begin{center}
\textbf{ А. Р. Миротин, А. А. Атвиновский}
\end{center}

\

Аннотация. Пусть функция  $f$ с вещественными  полюсами, образующими монотонную и ограниченную  последовательность, разлагается в ряд  Вольфа-Данжуа с положительными коэффициентами. Основной результат  заметки утверждает, что если  мы вычтем из функции $1/f$  ее  "линейную часть", то оставшаяся "дробная часть" этой функции тоже будет разлагаться в ряд  Вольфа-Данжуа (и ее полюсы тоже вещественны, а коэффициенты ряда отрицательны).
 Дано  приложение полученного результата к теории операторов.

Ключевые слова: ряд  Вольфа-Данжуа, замкнутый оператор, левый обратный оператор, функциональное исчисление.

\textbf{Введение}

Следуя \cite{Sib}, \textit{рядами Вольфа-Данжуа} мы будем называть ряды вида
$$
\sum\limits_{k=1}^{\infty}\frac{A_k}{z-\lambda_k},\eqno(\ast)
$$
где $A_k\in\mathbb{C}, \{A_k\}_{k\geq 1}\in l^1,$ $(\lambda_1, \lambda_2, \dots)$ --- ограниченная последовательность комплексных чисел.

 Как было отмечено в \cite{Sib}, ряды указанного вида интенсивно изучались в работах А.~Пуанкаре, Ж.~Вольфа, А.~Данжуа, Э.~Бореля,
 Т.~Карлемана, А.~Бёрлинга, Т.~А.~Леонтьевой \cite{L1}, \cite{L2}, прежде всего в связи с проблемами
 квазианалитичности и аналитического продолжения. Они также имеют приложения к теории рядов Дирихле и теории операторов \cite{Sib}.

Кроме того (см. статью  \cite{Sher} и приведенную там библиографию), различные свойства функций, допускающих представления вида $(\ast),$ изучались и использовались в ряде работ М.~Г.~Крейна, Г.~ Л.~ Гамбургера, Б.~ Я.~ Левина,
М.~ В.~ Келдыша и И.~ В.~ Островского, Л.~ де~ Бранжа,
Ю.~Ф.~ Коробейника, А.~ Боричева и М.~ Л.~ Содина, Л.~ С.~ Маергойза, В.~ Б.~ Шерстюкова по теории функций и гармоническому анализу, теории операторов и дифференциальных уравнений.    Отметим также работу Н.~ И.~ Ахиезера \cite{Ah}.

 Как известно, класс  функций Неванлинны $\mathcal{R}$ \cite{N} (см. также \cite{Don}, где эти функции называются функциями Пика) переходит в себя при преобразовании $f\mapsto -1/f.$
  Суть основного результата данной заметки состоит в уточнении этого
 свойства для некоторого подкласса класса  $\mathcal{R}.$    А именно,
 показано, что если функция  $f$ с вещественными полюсами, образующими монотонную и ограниченную  последовательность,  разлагается в ряд  Вольфа-Данжуа с положительными коэффициентами,  и если мы вычтем из  функции $1/f$ ее  "линейную часть", то оставшаяся "дробная часть" этой функции тоже будет разлагаться в ряд  Вольфа-Данжуа (и ее полюсы тоже вещественны, а коэффициенты ряда отрицательны). При этом оказалось, что условия положительности коэффициентов и вещественности полюсов нельзя отбросить (см.  замечания 1 и 2 ниже).
  Дано также приложение полученного результата к теории операторов.

   Статья опубликована в \cite{arx}.

\textbf{2 Основной результат}

\textbf{Теорема 1.} \textit{Пусть функция }$f$ \textit{представима в
виде}
$$
f(z)=\sum\limits_{k=1}^{\infty}\frac{c_k}{\lambda_k-z}, \eqno(1)
$$
\textit{где }$c_k>0, \sum_{k=i}^{\infty}c_k<\infty, \{\lambda_1,
\lambda_2, \dots\}$ \textit{--- монотонно возрастающая и ограниченная
последовательность действительных чисел.  Тогда}
$$
\frac{1}{f(z)}=\alpha+\beta z
-\sum_{n=1}^{\infty}\frac{b_n}{t_n-z},\eqno(2)
$$
\textit{где}
$$
\alpha=\frac{\sum_{k=1}^{\infty}c_k\lambda_k}{\left(\sum_{k=1}^{\infty}c_k\right)^2},\quad
\beta=-\frac{1} {\sum_{k=1}^{\infty}c_k},
$$
$t_n$ --- \textit{все нули функции} $f(z),$ $b_n=1/f'(t_n)>0\ (n=1,
2, \dots),$ $\sum_{n=1}^\infty b_n<\infty.$

Доказательство. Пусть $a:=\lambda_1, b:=\sup_{k} \lambda_k.$ Ясно, что  ряд (1) сходится локально равномерно и функция $f$ голоморфна на множестве $\mathbb{C}\setminus [a,b],$  в окрестности каждого из интервалов $(\lambda_{k-1}, \lambda_k)$  и  в бесконечности и имеет в бесконечности нуль первого порядка. Особыми точками функции $f$ в расширенной комплексной плоскости являются полюсы $\lambda_k$ и точка $b,$ предельная для полюсов.
Заметим, что все нули $t_n$ функции $f$  принадлежат отрезку $[a,b]$, так как ($z=x+iy$)
\[
 f(z)=\sum\limits_{k=1}^{\infty}\frac{c_k(\lambda_k-x)}{\left|\lambda_k-z\right|^2}+iy\sum\limits_{k=1}^{\infty}\frac{c_k}
 {\left|\lambda_k-z\right|^2},
\]
и кратность этих нулей равна единице, поскольку при $x\in \mathbb{R}\setminus(\{\lambda_k\}\cup \{b\})$
\[
f'(x)=\sum\limits_{k=1}^{\infty}\frac{c_k}{(\lambda_k-x)^2}>0.
\]
Из последнего неравенства следует также, что $f$ строго возрастает на любом интервале, содержащемся в $\mathbb{R}\setminus(\{\lambda_k\}\cup \{b\}).$  Поскольку $f(\lambda_{k-1}+0)=-\infty,$  $f(\lambda_{k}-0)=+\infty,$   множество  нулей функции $f$ на каждом интервале $(\lambda_{k-1}, \lambda_k)$ состоит ровно из одной точки. Следовательно, множество $\{t_n\}$ всех нулей функции $f$ счетно, и  мы можем считать последовательность $\{t_n\}$ строго возрастающей.

Покажем, что функция
$$
\varphi(z):=\frac{1}{f(z)}
$$
($\varphi(\lambda_k):=0,$ $\varphi(b):=0$)
имеет вид
\[
\varphi (z)=\alpha +\beta z-\int\limits _a^b \frac{d\tau(t)}{t-z},
\]
где  $\tau$ --- ограниченная неотрицательная регулярная борелевская мера, сосредоточенная  на $[a,b]$  (ср.  \cite{A}).
Действительно, так как
\[
\mathrm{Im}f(z)=y\sum\limits_{k=1}^{\infty}\frac{c_k}
 {\left|\lambda_k-z\right|^2},
\]
  то функция $-\varphi(z)$ принадлежит  классу Неванлинны $\mathcal{R}$ \cite{N} (см. также
  \cite[с.  217]{Don}), т.~е. голоморфна в верхней полуплоскости и $\mathrm{Im}(-\varphi(z))> 0$ при $y>0,$  а потому по теореме Неванлинны
\[
-\varphi(z)=\alpha _{1} +\beta _{1} z+\int\limits_{-\infty }^{+\infty }\left(\frac{1}{t-z}-\frac{t}{1+t^{2} } \right) d\tau(t),
\]
где $\tau$ --- неотрицательная регулярная борелевская мера, для которой сходится интеграл $\int _{-\infty }^{+\infty }(1+t^{2} )^{-1}d\tau(t)$,  $\alpha _{1}, \beta _{1}$ --- вещественные числа. Воспользуемся формулой обращения Стилтьеса-Перрона, которая в случае функции $-\varphi(z)$ при подходящей нормировке интегрирующей функции имеет вид (см., например, \cite[с. 521]{KN})
\[
\tau(s_{2} )-\tau(s_{1} )=-\frac{1}{\pi }\lim\limits_{\varepsilon \to 0} \int\limits _{s_1}^{s_2} \mathrm{Im}(\varphi(x+i\varepsilon ))dx
\]
(через $\tau(t)$ мы обозначаем функцию распределения  меры $\tau$).  Ввиду вещественности функции $\varphi(z)$ при $z=x<a$ и $z=x>b$ правая часть здесь равна нулю при $s_1<s_2<a$ и при $s_2>s_1>b$, а потому функция  $\tau(t)$  постоянна при $t<a$ и $t>b$ и, в частности, ограничена. Поэтому
\[
-\varphi(z)=\alpha_2+\beta _{1} z+\int\limits_a^b\frac{d\tau(t)}{t-z}  ,
\]
то есть
\[
\varphi(z)=\alpha +\beta z-\int\limits_a^b\frac{d\tau(t)}{t-z},
\]
где $\alpha, \beta$ --- вещественные числа.

Аналогично, если $t_{k-1}<s_1<s_2<t_k,$ формула обращения Стилтьеса-Перрона показывает, что  функция $\tau(t)$ постоянна на каждом интервале $(t_{k-1},t_k),$ что приводит к формуле (2), в которой $b_k>0$ --- скачок функции $\tau$  в точке $t_k$
(предельный переход под знаком интеграла в формуле Стилтьеса-Перрона возможен в силу непрерывности подинтегральной функции в комплексной окрестности отрезка $[s_1,s_2]$). Следовательно, $\sum_{k=1}^\infty b_k$ есть вариация функции  $\tau(t)$
на  отрезке $[a,b],$ а потому этот ряд сходится.

Наконец, из  равенства
$\lim_{x\to\infty}\sum_{n=1}^{\infty}b_n/(x-t_n)=0$
вытекает, что коэффициенты $\alpha$ и $\beta$ определяются следующим образом:
\[
\beta=\lim\limits_{x\to\infty}\frac{1}{xf(x)},\quad
\alpha=\lim\limits_{x\to\infty}\left( \frac{1}{f(x)}-\beta x\right).
\]
Тогда
$$
\beta=\lim\limits_{x\to \infty}\frac{1}{\sum_{k=1}^{\infty}\frac{xc_k}{\lambda_k -x}}=\lim\limits_{x\to \infty}\frac{1}
{\sum_{k=1}^{\infty}\frac{c_k}{\frac{\lambda_k}{x}-1}}=-\frac{1}{\sum_{k=1}^{\infty}c_k},
$$
и
$$
\alpha=\frac{1}{\sum_{k=1}^{\infty} c_k}\lim\limits_{x\to\infty}\frac{\sum_{k=1}^{\infty} c_k\left( 1+\frac{x}
{\lambda_k-x}\right)}{\sum_{k=1}^{\infty} \frac{c_k}{\lambda_k-x}}=
$$
$$
=\frac{1}{\sum_{k=1}^{\infty} c_k}\lim\limits_{x\to\infty}\frac{\sum_{k=1}^{\infty}\frac{c_k\lambda_k}{\lambda_k-x}}
{\sum_{k=1}^{\infty}\frac{c_k}{\lambda_k-x}}=\frac{1}{\sum_{k=1}^{\infty} c_k}\lim\limits_{x\to\infty}
\frac{\sum_{k=1}^{\infty}\frac{c_k\lambda_k}{\frac{\lambda_k}{x}-1}}{\sum_{k=1}^{\infty}\frac{c_k}{\frac{\lambda_k}{x}-1}}=
$$
$$
=\frac{\sum_{k=1}^{\infty} c_k\lambda_k}{\left(\sum_{k=1}^{\infty} c_k\right)^2},
$$
что и завершает доказательство теоремы.

\textbf{Замечание 1.} Требование положительности коэффициентов $c_k$ в теореме 1  существенно, так как в противном случае кратность нуля функции $f$ (а потому и порядок соответствующего полюса  функции  $1/f$) может быть больше единицы, что показывает следующий пример:
\[
\frac{1}{2z}-\frac{4}{z+1}+\frac{9}{2(z+2)}=\frac{(z-1)^2}{z(z+1)(z+2)}.
\]

\textbf{3. Приложения.}

Следствием  теоремы 1 является теорема об обращении функции $f$ вида (1) от замкнутого оператора в банаховом пространстве.
 Если $A$ --- замкнутый плотно определённый оператор в комплексном банаховом пространстве $X,$ спектр $\sigma(A)$
которого не пересекается с отрезком $[a,b],$ где $a:=\lambda_1, b:=\sup_{k} \lambda_k,$ то мы положим
$$
f(A):=\sum\limits_{k=1}^{\infty} c_kR(\lambda_k,A),
$$
где $R(\lambda_k,A)=(\lambda_k I-A)^{-1}$ --- значения резольвенты оператора $A.$ Это определение согласуется с голоморфным
функциональным исчислением Рисса-Данфорда замкнутых операторов в пространстве $X$ \cite{DS}, поскольку $f$ принадлежит пространству $\mathcal{F}(A)$
функций, голоморфных в некоторой (своей для каждой функции) окрестности множества $\sigma(A)$ и в бесконечности. Идущее ниже следствие обобщает результат из \cite{Obrresolvent}, который был получен с помощью
функционального исчисления, построенного в \cite{IZV1}---\cite{IZV2}. Континуальный аналог этого результата был установлен в \cite{Trudy}.

\textbf{Следствие 1. } \textit{Пусть  функция $f$ задана формулой (1),  $a:=\lambda_1, b:=\sup_{k} \lambda_k,$ и
$A$ --- замкнутый плотно определенный оператор в комплексном банаховом пространстве $X,$ спектр которого не пересекается
с отрезком $[a,b].$ Тогда левый обратный к оператору  $f(A)$ существует и имеет вид
 \[
f(A)^{-1}=\alpha I+\beta A-
\sum\limits_{n=1}^{\infty} b_n R(t_n,A),
 \]
где  $t_n\  (n=1,2,\dots)$ --- все нули функции   $f,$   а значения $b_n,$ $\alpha$  и
$\beta$  даются теоремой 1.}

Доказательство. Заметим, что функция $f$ принадлежит классу $R[a,b]$ (относительно последнего см.\cite[c. 525]{KN}, а также \cite{IZV1}). Поэтому по теореме 1  из \cite{IZV1} левый обратный к оператору  $f(A)$ существует и равен $\varphi(A)$ (как и выше,  $\varphi=1/f$). Осталось заметить, что в силу теоремы 1
 \[
\varphi(A)=\alpha I+\beta A-
\sum\limits_{n=1}^{\infty} b_n R(t_n,A).
 \]

В работе \cite{Obrresolvent} был рассмотрен частный случай следствия 1, в котором ряд  заменен конечной суммой. Там же была поставлена задача обобщения этого результата на случай комплексных полюсов (метод, использованный в \cite{Obrresolvent}, здесь неприменим, см. ниже замечание 2).  Нижеследующая теорема 2 решает эту задачу.

Рассмотрим рациональную функцию
$$
f(z)=\sum_{j=1}^n\frac{a_j}{\lambda_j-z},
$$
где $a_j>0$, а $\lambda_j$ --- произвольные попарно различные комплексные числа. Если $A$ --- замкнутый плотно определённый
оператор в комплексном банаховом пространстве $X$, спектр $\sigma(A)$ которого не пересекается с множеством
$\{\lambda_1,  \dots, \lambda_n\},$ то мы, как и выше, положим
$$
f(A)=\sum\limits_{j=1}^n a_jR(\lambda_j,A).
$$
 Это определение также согласуется с голоморфным
функциональным исчислением Рисса-Данфорда замкнутых операторов в пространстве $X,$ поскольку $f\in \mathcal{F}(A).$

Ниже мы получим условия левой обратимости оператора $f(A)$ и
вычислим соответствующий левый обратный.  Для формулировки основного результата заметим, что
рациональная функция $g=1/f$ имеет в бесконечности полюс первого
порядка. Следовательно, выделяя целую часть и разлагая дробную
часть на простейшие дроби мы можем её представить в виде:
$$
g(z)=\alpha+\beta z+\sum_{j=1}^m\sum_{k=1}^{m_j}\frac{c_{jk}}{(t_j-z)^k}, \eqno(3)
$$
где $t_j$ --- все нули функции $f$, $m_j$ --- кратность нуля $t_j$.

\textbf{Замечание 2.} В отличие от случая, когда $\lambda_j$ --- действительные числа, рассмотренного выше, нули
 функции $f$ могут быть кратными даже если все коэффициенты положительны. Например, так будет в случае
$$
f(z)=\frac{1}{z-\lambda}+\frac{1}{z-1}+\frac{1}{z},
$$
где $\lambda$ --- корень уравнения $\lambda^2-\lambda+1=0$.

\textbf{Теорема 2.} \textit{Пусть} $A$ \textit{--- замкнутый плотно определённый оператор в комплексном банаховом пространстве $X,$ спектр} $\sigma(A)$ \textit{которого не пересекается с выпуклой оболочкой}
 ${\rm conv}(\lambda_1,  \dots, \lambda_n)$ \textit{множества} $\{\lambda_1,  \dots, \lambda_n\}.$
 \textit{Тогда левый обратный к оператору} $f(A)$ \textit{существует и имеет вид}
$$
f(A)^{-1}=\alpha I+\beta A+\sum_{j=1}^m\sum_{k=1}^{m_j}c_{jk}R(t_j,A)^k,
$$
\textit{где} $t_j\  (j=1, \dots, m)$ \textit{--- все нули функции} $f$, $m_j$ \textit{--- кратность нуля} $t_j$ \textit{и}
$$
\alpha=\frac{\sum_{j=1}^n a_j\lambda_j}{\left(\sum_{j=1}^n a_j\right)^2}, \quad \beta=-\frac{1}{\sum_{j=1}^n a_j}.
$$

\textbf{Доказательство.} Значение коэффициентов $\alpha$ и $\beta$ выводятся из формулы (3) аналогично тому, как это было сделано в доказательстве теоремы 1.

Теперь покажем, что все корни уравнения $f(z)=0$ принадлежат ${\rm
conv}(\lambda_1,  \dots, \lambda_n).$ Если допустить
противное, то найдётся прямая на комплексной плоскости,
разделяющая ${\rm conv}(\lambda_1,  \dots, \lambda_n)$ и
некоторый корень $z_0$ этого уравнения. Следовательно, найдётся
прямая, разделяющая ${\rm conv}(\lambda_1-z_0,
\dots, \lambda_n-z_0)$ и $0$. Совершая поворот  $z\mapsto
e^{i\theta}z$ на подходящий угол, получаем, что прямая ${\rm
Re}f=a, a>0$ разделяет ${\rm conv}(e^{i\theta}(\lambda_1-z_0),
 \dots, e^{i\theta}(\lambda_n-z_0))$
и $0$. Ясно, что $\sum_{j=1}^n a_j/w_j=0,$
где $w_j=e^{i\theta}(\lambda_j-z_0).$

Дробно-линейное преобразование $\zeta=1/w$ переводит
прямую ${\rm Re}f=a$ в окружность, проходящую через $0$ и
содержащую внутри все точки $\zeta_j=1/w_j.$
Следовательно,
$$
\sum_{j=1}^n\frac{a_j}{w_j}=\sum_{j=1}^n a_j\zeta_j\ne 0.
$$
Мы получили противоречие (все $\zeta_j$  лежат по одну сторону от касательной  к окружности, проведенной  в точке $0$).

Из доказанного выше следует, что функция
$$
g(z)=\frac{1}{f(z)}=\alpha+\beta z+h(z),
$$
где
$$
h(z)=\sum_{j=1}^m\sum_{k=1}^{m_j}\frac{c_{jk}}{(t_j-z)^k},
$$
 голоморфна в окрестности спектра оператора $A,$ а потому функции $h$ принадлежат пространству
$\mathcal{F}(A).$ Из определения функционального исчисления Рисса-Данфорда сразу следует, что
$$
g(A)=\alpha I+\beta A+\sum_{j=1}^m\sum_{k=1}^{m_j}c_{jk}R(t_j,A)^k.
$$
Заметим, что оба слагаемых правой части очевидного равенства
$$
1=g(z)f(z)=(\alpha+\beta z)f(z)+h(z)f(z)
$$
принадлежат $\mathcal{F}(A).$ Следовательно, применяя функцию $(\alpha+\beta z)f(z)+h(z)f(z)$ к оператору $A$ и
воспользовавшись свойствами полиномиального исчисления  и голоморфного функционального исчисления Рисса-Данфорда \cite[VII.9]{DS}, будем иметь
$$
(\alpha I+\beta A)f(A)+h(A)f(A)=I.
$$
 Таким образом,
$f(A)^{-1}=g(A),$ что и требовалось доказать.

\textbf{Замечание 3.} Легко проверяемое равенство ($a, b\geq 0,$ $a+b=1$)
$$
aR(\lambda_1,A)+bR(\lambda_2,A)=R(\lambda_1,A)\left((a\lambda_2+b\lambda_1)-A\right)R(\lambda_2,A)
$$
показывает, что его левая часть  обратима слева тогда и только тогда, когда
число $a\lambda_2+b\lambda_1$ не принадлежат точечному спектру
$\sigma_p(A).$ Отсюда следует, что всевозможные линейные комбинации
с положительными коэффициентами двух значений резольвенты
оператора $A$ обратимы слева тогда и только тогда, когда $\sigma_p(A)$ не пересекается с выпуклой оболочкой ${\rm
conv}(\lambda_1, \lambda_2)$ множества
$\{\lambda_1,\lambda_2\}$, т.~ е. с отрезком с концами $\lambda_1$
и $\lambda_2$. Таким образом, условие
 ${\rm conv}(\lambda_1,  \dots, \lambda_n)\cap \sigma(A)=\emptyset$ в предыдущей теореме существенно.

\textbf{Замечание 4.}  Как в ситуации, описываемой следствием 1, так и в ситуации теоремы 2 задача $f(A)x=y$  является некорректной при неограниченном $A.$ Причиной этому служит член $\beta A$  в формуле для $f(A)^{-1},$ которая в обоих случаях имеет вид
$$
f(A)^{-1}=\alpha I+\beta A+g(A),
$$
где оператор $g(A)$  ограничен.
Из этой же формулы вытекает  следующий  подход  к регуляризации этих задач. Если оператор $A^{-1}$  ограничен и  $R_t$ ($0<t<t_0$) есть регуляризующее семейство
 задачи   $A^{-1}x=y$ (т.е. $R_t A^{-1}x\to x\ (t\to 0)$ при всех $x\in X,$ см., например, \cite{IVT}), то $R'_t=\alpha I+\beta R_t+g(A)$  есть регуляризующее семейство задачи   $f(A)x=y.$

\end{document}